\documentclass[12pt]{amsart}
\usepackage{graphicx}
\usepackage{amssymb}
\usepackage{amscd}
\usepackage{citesort}
\usepackage[dvips]{color}

\newtheorem{tw}{Theorem}[section]
\newtheorem{prop}[tw]{Proposition}
\newtheorem{lem}[tw]{Lemma}

\theoremstyle{remark}
\newtheorem{uw}[tw]{Remark}

\theoremstyle{definition}

\newcommand{\cal}[1]{\mathcal{#1}}
\newcommand{\bez}{\setminus}

\newcommand{\sig}{\sigma}
\newcommand{\eps}{\varepsilon}

\newcommand{\fal}[1]{\widetilde{#1}}
\newcommand{\dasz}[1]{\widehat{#1}}
\newcommand{\kre}[1]{\overline{#1}}
\newcommand{\gen}[1]{\langle #1 \rangle}
\newcommand{\map}[3]{#1\colon #2\to #3}
\newcommand{\field}[1]{\mathbb{#1}}
\newcommand{\zz}{\field{Z}}

\newcommand{\KL}[1]{ \bigl\{ #1 \bigr\}}
\newcommand{\st}{\;|\;}

\newcommand{\lst}[2]{{#1}_1,\dotsc,{#1}_{#2}}
\newcommand{\Mob}{M\"{o}bius strip}

\begin{document}

\numberwithin{equation}{section}

\title[Generating mapping class groups of nonorientable \ldots] {Generating mapping class groups of nonorientable surfaces with boundary}

\author{Micha\l\ Stukow}

\thanks{Supported by BW 5100-5-0205-6}
\address[]{
Institute of Mathematics, University of Gda\'nsk, Wita Stwosza 57, 80-952 Gda\'nsk, Poland }

\email{trojkat@math.univ.gda.pl}



\begin{abstract}
We obtain simple generating sets for various mapping class groups of a nonorientable surface
with punctures and/or boundary. We also compute the abelianizations of these mapping
class groups.
\end{abstract}

\maketitle%
\section{Introduction}%
Let $N_{g,s}^n$ be a smooth, nonorientable and compact surface of genus $g$ with $s$
boundary components and $n$ punctures. If $s$ and/or $n$ is zero then we omit it from the
notation. If we do not want to emphasise the numbers $g,s,n$, we simply write $N$ for a
surface $N_{g,s}^n$. Recall that $N_{g}$ is a connected sum of $g$ projective planes and
$N_{g,s}^n$ is obtained from $N_g$ by removing $s$ open discs and specifying the set
$\Sigma=\{\lst{z}{n}\}$ of $n$ distinguished points in the interior of $N$.

Let ${\textrm{Diff}}(N)$ be the group of all diffeomorphisms $\map{h}{N}{N}$ such that
$h$ is the identity on each boundary component and $h(\Sigma)=\Sigma$. By ${\cal{M}}(N)$
we denote the quotient group of ${\textrm{Diff}}(N)$ by the subgroup consisting of maps
isotopic to the identity, where we assume that isotopies fix $\Sigma$ and are the
identity on each boundary component. ${\cal{M}}(N)$ is called the \emph{mapping class
group} of $N$. The mapping class group of an orientable surface is defined analogously,
but we consider only orientation preserving maps.

By abuse of notation we will use the same letter for a map and its isotopy class and we
will use the functional notation for the composition of diffeomorphisms.

For any $1\leq k\leq n$, let ${\cal{PM}}^k(N)$ be the subgroup of ${\cal{M}}(N)$
consisting of elements which fix $\Sigma$ pointwise and preserve the local orientation
around the punctures $\{z_1\ldots,z_k\}$. For $k=0$, we obtain the so--called \emph{pure
mapping class group} ${\cal{PM}}(N)$.



The main reason for introducing the groups ${\cal{PM}}^k(N)$ with $k\geq 1$ is as
follows. The usual way to study the mapping class group of a surface with boundary is via
the homomorphism
\[\map{i_*}{{\cal{PM}}(N_{g,s}^n)}{{\cal{PM}}(N_{g}^{n+s})},\]
induced by the inclusion $\map{i}{N_{g,s}^n}{N_{g}^{n+s}}$, where $N_{g}^{n+s}$ is the
surface obtained by gluing a disk with a puncture to each boundary component of
$N_{g,s}^{n}$. However, in the nonorientable case, the homomorphism $i_*$ is not an
epimorphism. To be more precise, its image consists of maps which preserve the local
orientation around the $s$ punctures coming from the boundary components of $N_{g,s}^n$.
Hence the group ${\cal{PM}}^k(N)$ occurs naturally as an image of $i_*$.

%
%
%
\subsection{Background}
One of the oldest problems concerning mapping class groups is to find a simple generating
set or a generating set with some additional properties. In the case of an orientable
surface there are many results in this direction -- see for example \cite{Bir-Punct,ParLab,Lick1,Lick2,McCarPap1,Kork-twogen,Hump,Bre-Farb,Stukow_HypGen,MacMod,Stukow_inv,Wajn} and references there.

On the other hand, the nonorientable case has not been studied much. The first
significant result is due to Lickorish \cite{Lick3}, who proved that the mapping class
group ${\cal{M}}(N_{g})$ is generated by Dehn twists and a so--called ``crosscap slide"
(or a ``Y--homeomorphism"). Later Chillingworth \cite{Chil} found a finite generating set for the group ${\cal{M}}(N_{g})$, which was extended by Korkmaz \cite{Kork-non} to the case of groups ${\cal{M}}(N_{g}^n)$ and ${\cal{PM}}(N_{g}^n)$.
It is also known that the group
${\cal{M}}(N_{g}^n)$ is generated by involutions \cite{Szep1,Szep2}.

Another natural question is to compute the (co)homology groups of mapping class groups.
As above, there are many results concerning the orientable case -- a good reference is a
survey article \cite{Kork-hom}. In the nonorientable case, Korkmaz
\cite{Kork-non1,Kork-non} computed the first integral homology group of ${\cal{M}}(N_{g}^n)$, and
under additional assumption $g\geq 7$, of ${\cal{PM}}(N_{g}^n)$.
\subsection{Main results}
The main goal of this paper is to extend some of the above results to the case of mapping
class groups of nonorientable surfaces with boundary. More precisely, for every $g\geq 3$
we obtain finite generating sets for the groups ${\cal{PM}}^k(N_{g,s}^n)$ and
${\cal{M}}(N_{g,s}^n)$ -- cf Theorems \ref{tw:gen:pure:bd} and \ref{tw:gen:notpure}. Then
using these generating sets we compute their first integral homology groups
(abelianizations)  -- cf Theorems \ref{tw:hom:pure} and \ref{tw:hom:notpure}.

The reason for the assumption $g\geq 3$ is the exceptional (and nontrivial) nature of the cases $g=1$ and $g=2$. Most of our analysis make no sense in these cases, hence we leave them for future consideration.

Let us point out that although, using the results of \cite{BirChil1}, one can prove that the mapping class group of a closed
nonorientable surface is an infinite index subgroup of the mapping class group of an
orientable surface, we do not see any method to deduce our results from the orientable case.

It is worth mentioning that surfaces with boundary occur in a natural way when one
considers the stabilisers in the mapping class group of sets of circles. Such a situation
is very common when dealing with complexes of curves on a surface. In particular we
believe that our work will be an important step toward finding a presentation for the
mapping class group of a nonorientable surface \cite{Szep_curv}.

\section{Preliminaries}
By a \emph{circle} on $N$ we mean an oriented simple closed curve on $N\bez \Sigma$, which is
disjoint from the boundary of $N$. Usually we identify a circle with its image. Moreover, as in the case of diffeomorphisms, we will use the same letter for a circle and its isotopy class. According to whether a regular neighbourhood of
a circle is an annulus or a \Mob, we call the circle \emph{two--sided} or \emph{one--sided} respectively.
We say that a circle is \emph{generic} if it bounds neither a disk with less than $2$ punctures nor a \Mob\ disjoint from $\Sigma$.

Let $a$ be a two--sided circle. By definition, a regular neighbourhood of $a$ is an
annulus, so if we fix one of its two possible orientations, we can define the
\emph{right Dehn twist} $t_a$ about $a$ in the usual way. We emphasise that since we are
dealing with nonorientable surfaces, there is no canonical way to choose the orientation
of $S_a$. Therefore by a twist about $a$ we always mean one of two possible twists about
$a$ (the second one is then its inverse). By a \emph{boundary twist} we mean a twist
about a circle isotopic to a boundary component. It is known that if $a$ is not
generic then the Dehn twist $t_a$ is trivial. In particular, a Dehn twist about the
boundary of a \Mob\ is trivial -- see Theorem 3.4 of \cite{Epstein}.

Other important examples of diffeomorphisms of a nonorientable surface are the
\emph{crosscap slide} and the \emph{puncture slide}. They are defined as a slide of a
crosscap and of a puncture respectively, along a loop. The general convention is that one
considers only crosscap slides along one--sided simple loops (in such a form they were
introduced by Lickorish \cite{Lick3}), for precise definitions and properties see
\cite{Kork-non}.

The following two propositions follow immediately from the above definitions.

\begin{prop}\label{pre:prop:twist:cong}
Let $N_a$ be an oriented regular neighbourhood of a two--sided circle $a$ in a surface
$N$, and let $\map{f}{N}{N}$ be any diffeomorphism. Then $ft_af^{-1}=t_{f(a)}$, where the
orientation of a regular neighbourhood of $f(a)$ is induced by the orientation of
$f(N_a)$.\qed
\end{prop}

\begin{prop}\label{pre:prop:slide:cong}
Let $v$ be a slide of a puncture $z$ along a simple closed loop $\alpha$ on a surface
$N$, and let $\map{f}{N}{N}$ be any diffeomorphism. Then $fvf^{-1}$ is a slide of the
puncture $f(z)$ along the loop $f(\alpha)$. \qed
\end{prop}

The next proposition, which provides a relationship between puncture slides and twists is
proved in Section 6.1 of \cite{Ivanov1}.
\begin{prop}\label{pre:prop:slide:twist}
Let $\alpha$ be a two--sided simple loop on a surface $N$, based at the puncture $z$.
Define also $a$ and $b$ to be the boundary circles of a regular neighbourhood $N_\alpha$
of $\alpha$ such that the orientations of $\alpha$ and of $N_\alpha$ are as in Figure
\ref{01_pre} (we indicate the orientation of $N_\alpha$ by choosing the direction of a
right twist about $a$).
\begin{figure}[h]
\includegraphics{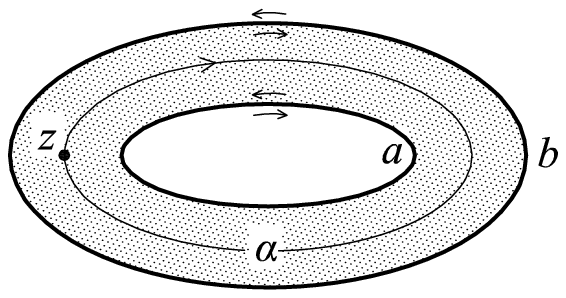}
\caption{Two--sided loop $\alpha$ and its regular neighbourhood.}\label{01_pre}
\end{figure}
Then $t_at_b^{-1}$ is the slide of $z$ along~$\alpha$. \qed
\end{prop}

Finally, let us recall the so--called \emph{lantern relation}, which will be our main
tool in studying properties of mapping class groups. The proof can be found in Section 4
of \cite{John1}.
\begin{prop}
Let $S$ be a sphere with four holes embedded in a surface $N\bez \Sigma$ and let
$a_0,a_1,a_2,a_3$ be the boundary circles of $S$. Define also $a_{1,2},a_{2,3},a_{1,3}$
as in Figure \ref{01_2_pre} and assume that the orientations of regular neighbourhoods of
these seven circles are induced from the orientation of $S$.
\begin{figure}[h]
\includegraphics{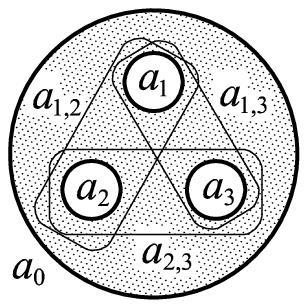}
\caption{Circles of the lantern relation.}\label{01_2_pre}
\end{figure}
Then \[t_{a_0}t_{a_1}t_{a_2}t_{a_3}=t_{a_{1,2}}t_{a_{2,3}}t_{a_{1,3}}.\] \qed
\end{prop}

\section{Generators for the group ${\cal{PM}}(N_g^n)$} 
\begin{figure}[h]
\includegraphics{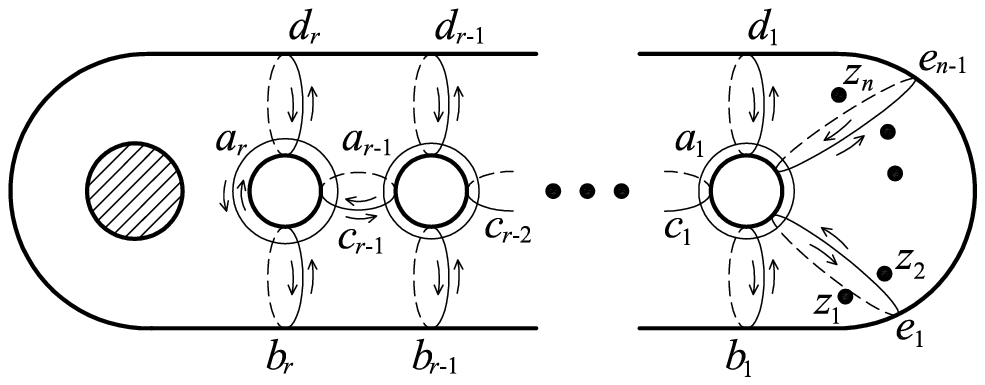}
\caption{Circles ${\cal{C}}$ for $g=2r+1$.}\label{fig:02_GenRed}
\end{figure}

\begin{figure}[h]
\includegraphics{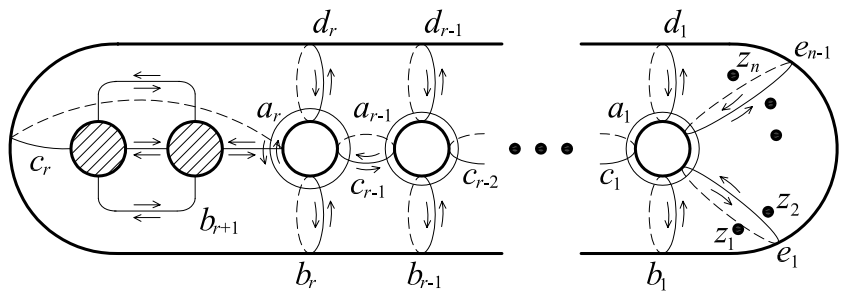}
\caption{Circles ${\cal{C}}$ for $g=2r+2$.}\label{fig:03_GenRed}
\end{figure}

Let us recall a known generating
set for the group ${\cal{PM}}(N_g^n)$ with $g\geq 3$. Represent the surface $N=N_{g}^n$
as a connected sum of an orientable surface and one or two projective planes (one for $g$
odd and two for $g$ even). Figures \ref{fig:02_GenRed} and \ref{fig:03_GenRed} show this
model of $N$ -- in these figures the shaded disks represent crosscaps, hence their
interiors are to be removed and then the antipodal points on each boundary component are
to be identified.

Let ${\cal{C}}$ be the set of circles indicated in Figure \ref{fig:02_GenRed} for $g=2r+1$,
and in Figure \ref{fig:03_GenRed} for $g=2r+2$.
Hence
\[ {\cal{C}}=\{\lst{a}{r},\lst{b}{r},\lst{c}{r-1},\lst{d}{r},\lst{e}{n-1}\},\]
for $g=2r+1$, and
\[ {\cal{C}}=\{\lst{a}{r},\lst{b}{r+1},\lst{c}{r},\lst{d}{r},\lst{e}{n-1}\},\]
for $g=2r+2$. The figures also indicate our chosen orientations of local neighbourhoods
of circles in ${\cal{C}}$, the orientation is such that the arrow points to the right if we
approach the circle. Therefore by a twist about one of the circles in ${\cal{C}}$ we will
always mean the twist determined by this particular choice of orientation (recall that
the general rule is that we consider \emph{right} Dehn twists, i.e. if we approach the
circle of twisting we turn to the right). In what follows we will often indicate the
orientation of a regular neighbourhood of a circle by drawing the direction of a twist.

\begin{figure}[h]
\includegraphics{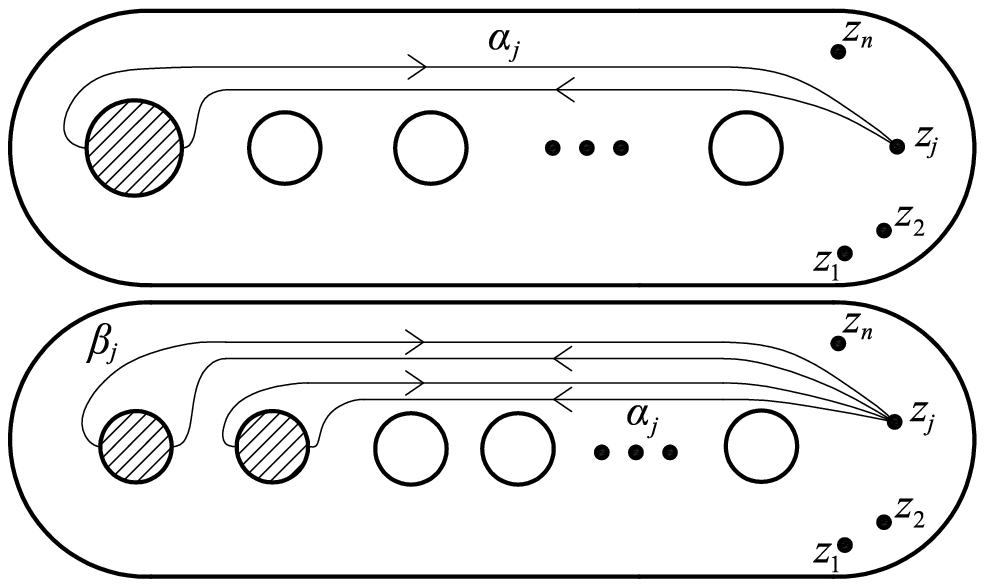}
\caption{Loops $\lst{\alpha}{n}$ and $\lst{\beta}{n}$.}\label{fig:04_GenRed}
\end{figure}

Let $v_j$ be a slide of a puncture $z_j$ along the loop $\alpha_j$ for $j=1,\ldots,n$ as
in Figure \ref{fig:04_GenRed}.
If $g=2r+2$, let $\lst{w}{n}$ be puncture slides along $\lst{\beta}{n}$ -- cf Figure
\ref{fig:04_GenRed}.

\begin{figure}[h]
\includegraphics{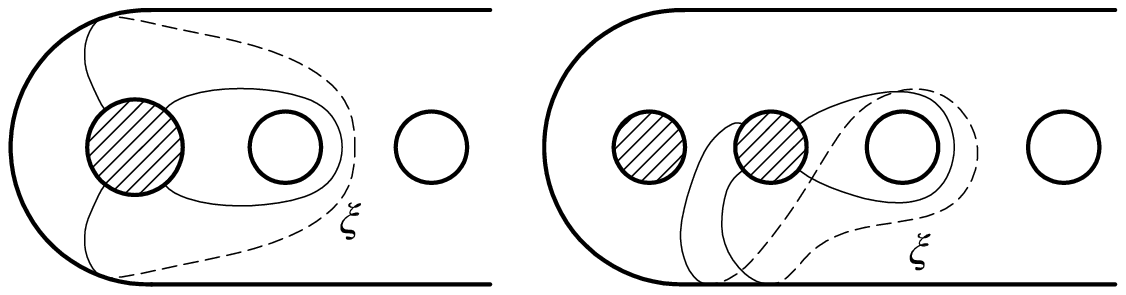}
\caption{Circle $\xi$ for $g=2r+1$ and $g=2r+2$.}\label{fig:05_GenRed}
\end{figure}

Define also $y$ to be a crosscap slide such that $y^2$ is a twist about the circle $\xi$
indicated in Figure \ref{fig:05_GenRed}. To be more descriptive, let $N_1$ be the
connected component of $N\bez \xi$ diffeomorphic to a Klein bottle with one boundary
component. Then $N_1$ is diffeomorphic to a disk with two crosscaps and we can define $y$
to be a slide of one of these crosscaps along the core of the second one. It turns out
that in what follows, the ambiguity in the definition of $y$ is inessential.
\begin{tw}[Theorem 4.13 of \cite{Kork-non}]\label{tw:kork:gen} Let $g\geq 3$. Then the mapping class group
${{\cal{PM}}(N_g^n)}$ is generated by
\begin{itemize}
 \item $\{t_l,v_j,y\st{l\in\cal{C}},1\leq j\leq n\}$ if $g$ is odd and
 \item $\{t_l,v_j,w_j,y\st{l\in\cal{C}},1\leq j\leq n\}$ if $g$ is even.
\end{itemize}\qed
\end{tw}
Now let us simplify the above generating set for $g$ even (we will replace all the $w_j$'s by a
single twist).

\begin{figure}[h]
\includegraphics{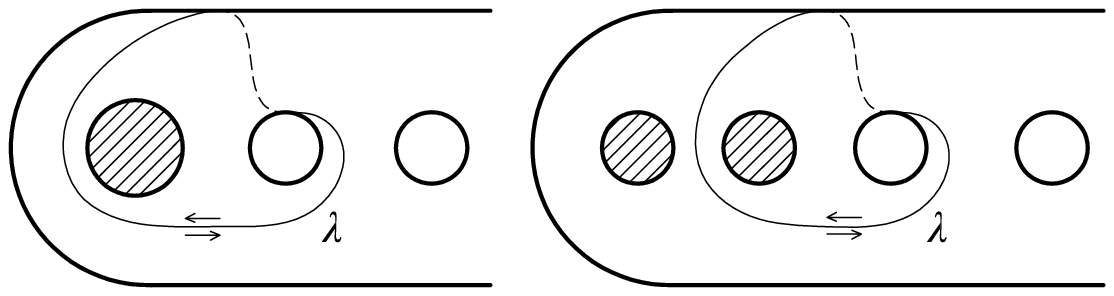}
\caption{Circle $\lambda$ for $g=2r+1$ and $g=2r+2$.}\label{fig:06_GenRed}
\end{figure}

Let $\lambda$ be the circle indicated in Figure \ref{fig:06_GenRed}.
\begin{tw} \label{tw:gen:pure:cl:even}
Let $g=2r+2\geq 4$. Then the mapping class group ${\cal{PM}}(N_g^n)$ is generated by
 \[\{t_l,v_j,y,t_\lambda\st{l\in{\cal{C}}},1\leq j\leq n\}.\]
\end{tw}
\begin{proof}
Let $G$ denote the group generated by the above elements. By Theorem \ref{tw:kork:gen},
it is enough to prove that $w_j\in G$ for $j=1,\ldots,n$. Our first claim is that if
$\delta_j$ is a circle as in Figure \ref{fig:07_GenRed}(ii) then
\begin{equation}
(t_{b_{r+1}}^{-1}(w_j)t_{b_{r+1}})v_j^{-1}=t_{\delta_j}^{-1}. \label{eq:w_i_kw}
\end{equation}
For notational convenience, for any simple loop $\alpha$ based at $z_j$, let $p(\alpha)$
be the slide of $z_j$ along $\alpha$. By Proposition \ref{pre:prop:slide:cong}, the
left--hand side of \eqref{eq:w_i_kw} can be rewritten as follows (observe that we compose
loops from left to right).
\[(t_{b_{r+1}}^{-1}(w_j)t_{b_{r+1}})v_j^{-1}=p(t_{b_{r+1}}^{-1}(\beta_j))p(\alpha^{-1})=
p(\alpha^{-1}t_{b_{r+1}}^{-1}(\beta_j))\]
 It is not hard to check that Figure \ref{fig:07_GenRed}(i) shows the
loop $t_{b_{r+1}}^{-1}(\beta_j)$. Hence $\alpha^{-1}t_{b_{r+1}}^{-1}(\beta_j)$ is the
loop shown in Figure \ref{fig:07_GenRed}(ii). By Proposition \ref{pre:prop:slide:twist},
the slide along this loop is equal to $t_{\delta_j}^{-1}$ which completes the proof of
\eqref{eq:w_i_kw}.
%

\begin{figure}[h]
\includegraphics{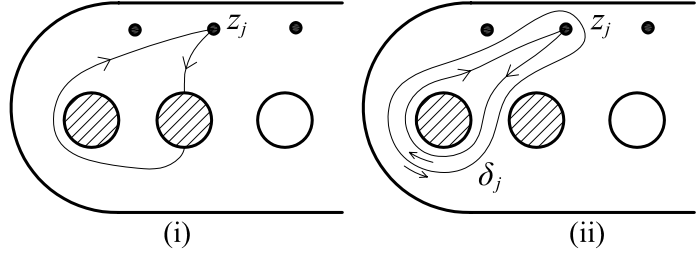}
\caption{Loops $t_{b_{r+1}}^{-1}(\beta_j)$, $\alpha_j^{-1} t_{b_{r+1}}^{-1}(\beta_j)$ and
circle $\delta_j$.}\label{fig:07_GenRed}
\end{figure}
Therefore, by equation \eqref{eq:w_i_kw}, it is enough to prove that $t_{\delta_j}\in G$.
Before we do that we need two lemmas.
\begin{lem}\label{lem:e:nu}
Let $e_{i,j}$ and $\nu_{i,m}$ be the circles shown in Figure \ref{fig:08_GenRed} for
$i=1,\ldots,r$, $j=0,\ldots,n$ and $m=1,\ldots,n$. Then the twists $t_{e_{i,j}}$ and
$t_{\nu_{i,m}}$ are in the group generated by $\{t_l\st l\in{\cal{C}}\}$.
\begin{figure}[h]
\includegraphics{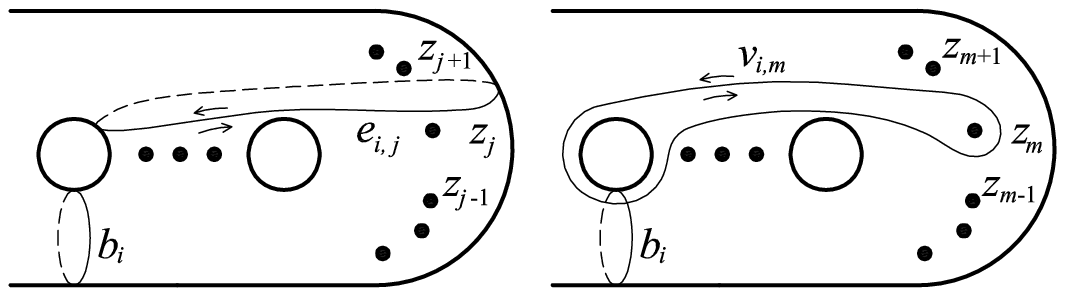}
\caption{Circles $e_{i,j}$ and $\nu_{i,m}$.}\label{fig:08_GenRed}
\end{figure}
\end{lem}
\begin{proof}
Let $H=\gen{t_l\st l\in{\cal{C}}}$. It is straightforward to check that
\[e_{i,j}
=t_{a_i}t_{c_{i-1}}t_{a_{i-1}}t_{b_i}^{-1}t_{a_i}^{-1}t_{c_{i-1}}^{-1}t_{b_{i-1}}t_{a_{i-1}}^{-1}(e_{i-1,j}).
\]
Moreover, $e_{0,0}=b_1$, $e_{0,n}=d_1$ and $e_{0,i}=e_i$ for $i=1,\ldots,n-1$. Therefore,
by induction on $i$, $t_{e_{i,j}}\in H$ (cf Proposition \ref{pre:prop:twist:cong}).

The rest of the proof follows, by Proposition \ref{pre:prop:twist:cong}, from the
relation
\[\nu_{i,m}=t_{e_{i,m-1}}t_{a_i}(e_{i,m}).\]
\end{proof}
\begin{uw}\label{rem:e:nu}
For further reference, observe that neither in the definition of the circles $e_{i,j}$ and
$\nu_{i,m}$ nor in the proof of the above lemma, we used the assumption that $g$ is even.
\end{uw}
\begin{lem}\label{lem:kap:tau}
Let $\tau$ and $\tau_j$ be the circles indicated in Figure \ref{fig:10_GenRed} for
$j=1,\ldots,n$. Then $t_{\tau},t_{\tau_j}\in G$.
\end{lem}
\begin{proof}
By Proposition \ref{pre:prop:twist:cong} and Lemma \ref{lem:e:nu}, the assertion follows
from the relations
\begin{align*}
\tau&=t_{d_r}t_{a_r}^{2}t_{b_r}t_{d_r}(\lambda),\\
\tau_j&=t_{e_{r,j-1}}t_{e_{r,j}}^{-1}(\tau).
\end{align*} 
\end{proof}
Now let us come back to the proof of Theorem \ref{tw:gen:pure:cl:even}. As was observed,
it is enough to prove that $t_{\delta_j}\in G$. Observe that the seven circles indicated in
Figure \ref{fig:10_GenRed} form a configuration of the lantern relation (the four circles
in Figure \ref{fig:10_GenRed}(i) bound a sphere with four holes). Therefore we have the
relation
\[t_{a_r}t_{\tau_j}=t_{\delta_j}t_{\tau}t_{\nu_{r,j}}.\]
By Lemmas \ref{lem:e:nu} and \ref{lem:kap:tau}, this implies that $t_{\delta_j}\in G$.
\begin{figure}[h]
\includegraphics{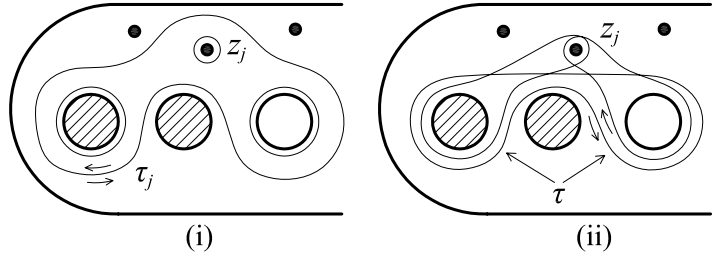}
\caption{Circles of the lantern relation
$t_{a_r}t_{\tau_j}=t_{\delta_j}t_{\tau}t_{\nu_{r,j}}$.}\label{fig:10_GenRed}
\end{figure}
\end{proof}

\section{Generators for the group ${\cal{PM}}^k(N_g^n)$} The following proposition can be
found in any book on combinatorial group theory -- see for example Chapter 9 of \cite{Johnson}.
\begin{prop}\label{prop:Joh}
Let $X$ be a generating set for a group $G$ and let $U$ be a left transversal for a subgroup
$H$ (i.e. $U$ is a set of representatives of left cosets of $H$) . Then $H$ is generated by the set
\[\{ux\kre{ux}^{-1}\, :\, u\in U, x\in X, ux\not\in U\},\]
where $\kre{g}=gH\cap U$ for $g\in G$. \qed\end{prop}

\begin{figure}[h]
\includegraphics{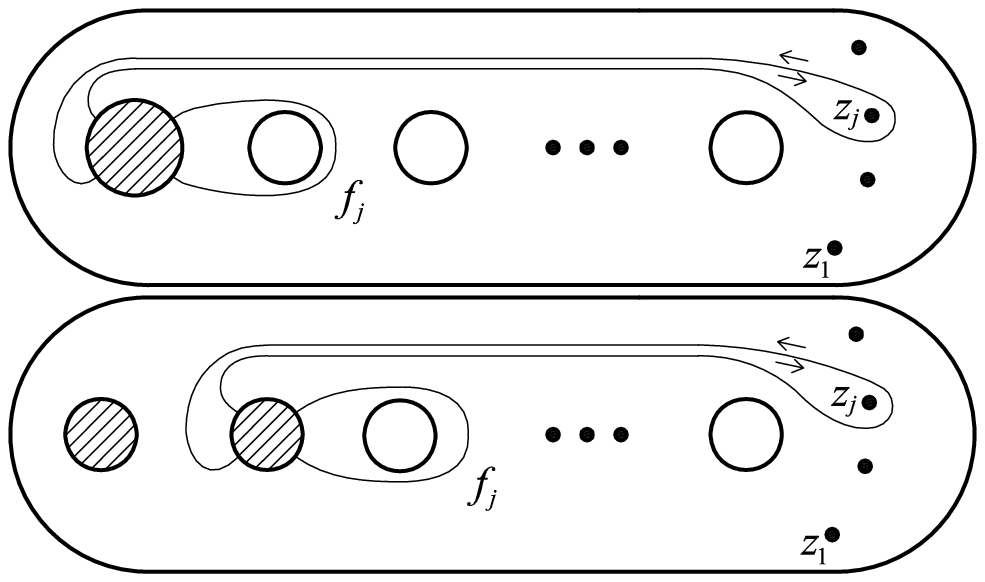}
\caption{Circles $\lst{f}{n}$.}\label{fig:11_GenCl}
\end{figure}

Let $\lst{f}{n}$ be the circles indicated in Figure \ref{fig:11_GenCl}.
\begin{tw}\label{tw:gen:pure:cl}
Let $g \geq 3$ and $0\leq k\leq n$. Then the mapping class group ${\cal{PM}}^k(N_g^n)$ is
generated by
\begin{itemize}
 \item $\{t_l,t_{f_1},\ldots,t_{f_k},v_{k+1},\ldots,v_n,y\st{l\in{\cal{C}}}\}$ if $g$ is
 odd and
 \item $\{t_l,t_{f_1},\ldots,t_{f_k},v_{k+1},\ldots,v_n,y,t_{\lambda}\st{l\in{\cal{C}}}\}$ if $g$ is even.
\end{itemize}
\end{tw}
\begin{proof}
The proof will be by induction on $k$. For $k=0$ the theorem follows from Theorems
\ref{tw:kork:gen} and \ref{tw:gen:pure:cl:even}. Suppose that the theorem is true for
$k-1$, i.e. the group $G={\cal{PM}}^{k-1}(N_g^n)$ is generated by the set
\begin{itemize}
 \item $\{t_l,t_{f_1},\ldots,t_{f_{k-1}},v_{k},\ldots,v_n,y\st{l\in{\cal{C}}}\}$ if $g$ is
 odd and
 \item $\{t_l,t_{f_1},\ldots,t_{f_{k-1}},v_{k},\ldots,v_n,y,t_{\lambda}\st{l\in{\cal{C}}}\}$ if $g$ is even.
\end{itemize}
If $H={\cal{PM}}^{k}(N_g^n)$, then $H$ is of index two in $G$, hence as a transversal for
$H$ we can take $U=\{1,v_k\}$. By Proposition \ref{prop:Joh}, $H$ is generated by
\begin{itemize}
 \item
 $\{t_l,v_kt_lv_k^{-1},t_{f_1},\ldots,t_{f_{k-1}}, v_kt_{f_1}v_k^{-1},\ldots,v_kt_{f_{k-1}}v_k^{-1},\\
 v_{k+1},\ldots,v_n,v_k^2, v_kv_{k+1}v_k^{-1},\ldots,v_kv_nv_k^{-1},
 y,v_kyv_k^{-1}\st{l\in{\cal{C}}}\}$\\
 if $g$ is odd and
 \item  $\{t_l,v_kt_lv_k^{-1},t_{f_1},\ldots,t_{f_{k-1}}, v_kt_{f_1}v_k^{-1},\ldots,v_kt_{f_{k-1}}v_k^{-1},\\
 v_{k+1},\ldots,v_n,v_k^2, v_kv_{k+1}v_k^{-1},\ldots,v_kv_nv_k^{-1},
 y,v_kyv_k^{-1},\\t_\lambda,v_k t_\lambda v_k^{-1}\st{l\in{\cal{C}}}\}$
 if $g$ is even.
\end{itemize}
Let $K\leq G$ be the group generated by
\begin{itemize}
 \item $\{t_l,t_{f_1},\ldots,t_{f_k},v_{k+1},\ldots,v_n,y\st{l\in{\cal{C}}}\}$ if $g$ is
 odd and
 \item $\{t_l,t_{f_1},\ldots,t_{f_k},v_{k+1},\ldots,v_n,y,t_{\lambda}\st{l\in{\cal{C}}}\}$ if $g$ is even.
\end{itemize}
Since each generator of $K$ preserves the local orientation around each of $\lst{z}{k}$,
we have $K\leq H$. To complete the proof it is enough to show that $H\leq K$, i.e. that
each generator of $H$ is in $K$. Since $v_k$ commutes with
$t_{a_i},t_{b_i},t_{c_j},t_{e_m},t_{\lambda}$ for $i=1,\ldots,r$, $j=1,\ldots,r-1$,
$m=1,\ldots,k-1$, it is enough to prove that
\begin{itemize}
 \item
 $v_kt_{d_1}v_k^{-1},\ldots,v_kt_{d_r}v_k^{-1},\\
 v_kt_{e_{k}}v_k^{-1},\ldots,v_kt_{e_n}v_k^{-1},\\
 v_kt_{f_1}v_k^{-1},\ldots,v_kt_{f_{k-1}}v_k^{-1},\\
v_kv_{k+1}v_k^{-1},\ldots,v_kv_{n}v_k^{-1},\\ v_k^2, v_kyv_k^{-1}\in K\quad \text{if $g$
is odd and}$
 \item
 $v_kt_{d_1}v_k^{-1},\ldots,v_kt_{d_r}v_k^{-1},\\
 v_kt_{e_{k}}v_k^{-1},\ldots,v_kt_{e_n}v_k^{-1},\\
 v_kt_{b_{r+1}}v_k^{-1},v_kt_{c_r}v_k^{-1}\\
 v_kt_{f_1}v_k^{-1},\ldots,v_kt_{f_{k-1}}v_k^{-1},\\
v_kv_{k+1}v_k^{-1},\ldots,v_kv_{n}v_k^{-1}, \\
v_k^2, v_kyv_k^{-1}\in K \quad \text{if $g$ is even.}$
\end{itemize}
The rest of this section is devoted to the proof of the above statements.
\subsection{Case of $v_kt_{d_i}v_k^{-1}$}
\begin{lem}\label{lem:mu}
Let the circles $\mu_{i,j}$ for $i=1,\ldots,r$, $j=1,\ldots,k$, be as in Figure
\ref{fig:12_GenCl}. Then $t_{\mu_{i,j}}\in K$.
\begin{figure}[h]
\includegraphics{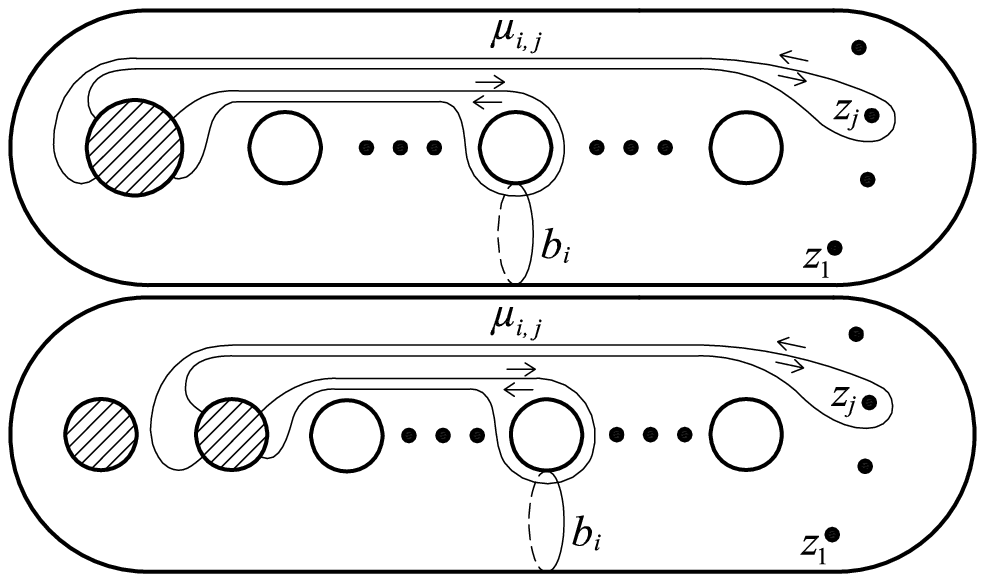}
\caption{Circles $\mu_{i,j}$.}\label{fig:12_GenCl}
\end{figure}
\end{lem}
\begin{proof}
One can check that
\[\mu_{i,j}=
t_{c_{i-1}}t_{a_{i}}t_{a_{i-1}}t_{c_{i-1}}(\mu_{i-1,j})\quad \text{for $i=2,\ldots,r$,
$j=1,\ldots,k$.} \] Moreover, $\mu_{r,j}=f_j$. Therefore, by Proposition
\ref{pre:prop:twist:cong}, the lemma follows by descending induction on $i$.
\end{proof}
By Lemma \ref{lem:e:nu} and Remark \ref{rem:e:nu}, $t_{\nu_{i,j}}\in
K$, for $i=1,\ldots,r$, $j=1,\ldots,n$. Therefore by Lemma \ref{lem:mu}, Proposition
\ref{pre:prop:twist:cong}, and by the relation
\[v_k(d_i)=t_{\mu_{i,k}}^{-1}t_{d_i}(\nu_{i,k}),
\] we obtain $v_kt_{d_i}v_k^{-1}\in K$, for $i=1,\ldots,r$.
\subsection{Case of $v_kt_{e_j}v_k^{-1}$} As above, by the relation
\[v_k(e_j)=t_{\mu_{1,k}}^{-1}t_{e_j}(\nu_{1,k}),\]
we obtain $v_kt_{e_j}v_k^{-1}\in K$, for $j=k,\ldots,n$.
\subsection{Case of $v_kt_{b_{r+1}}v_k^{-1}$, $g$ -- even}
Let $\rho=t_{d_r}t_{a_r}(c_r)$ and $\rho_k=t_{e_{r,k}}^{-1}t_{e_{r,k-1}}(\rho)$ (cf
Figure \ref{fig:13_GenCl}). By Lemma \ref{lem:e:nu} and Remark \ref{rem:e:nu}, $t_\rho,
t_{\rho_k}\in K$.
\begin{figure}[h]
\includegraphics{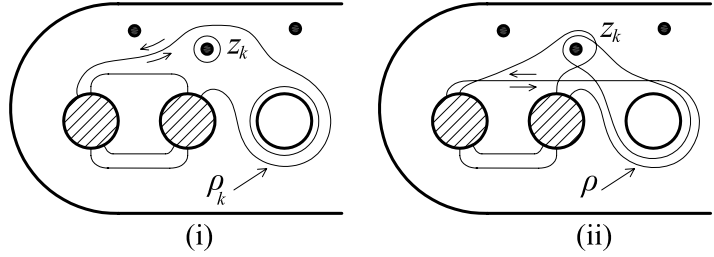}
\caption{Circles of the lantern relation
$t_{a_r}t_{\rho_k}t_{b_{r+1}}=t_{v_k({b_{r+1}})}t_\rho
t_{\nu_{r,k}}$.}\label{fig:13_GenCl}
\end{figure}
Now the seven circles indicated in Figure \ref{fig:13_GenCl} form a configuration of the
lantern relation (the four circles in Figure \ref{fig:13_GenCl}(i) bound a sphere with four
holes). Moreover, one of the circles in Figure \ref{fig:13_GenCl}(ii) is
$v_k({b_{r+1}})$. Therefore we have the relation
\[t_{a_r}t_{\rho_k}t_{b_{r+1}}=t_{v_k({b_{r+1}})}t_\rho t_{\nu_{r,k}}.\]
This proves that $t_{v_k({b_{r+1}})}=v_kt_{b_{r+1}}v_k^{-1}\in K$.
\subsection{Case of $v_k^2$} \label{sec:v_k}
First observe that if $\lambda$ is as in Figure \ref{fig:06_GenRed}, then $t_\lambda\in
K$. For $g$ even this follows from the definition of $K$, and for $g$ odd we have the
relation $\lambda=t_{b_r}^{-1}(a_r)$. Therefore if $\omega=t_{d_r}(\lambda)$ and
$\omega_k=t_{e_{r,k-1}}t_{e_{r,k}}^{-1}(\omega)$ (cf Figure \ref{fig:14_GenCl}) then
$t_\omega,t_{\omega_k}\in K$.
\begin{figure}[h]
\includegraphics{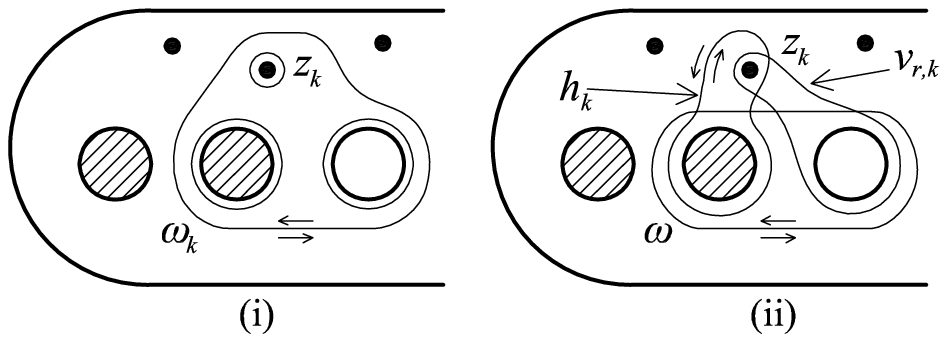}
\caption{Circles of the lantern relation $t_{a_r}t_{\omega_k}=t_{h_k}t_\omega
t_{\nu_{r,k}}$.}\label{fig:14_GenCl}
\end{figure}
Now the seven circles indicated in Figure~\ref{fig:14_GenCl} form a configuration of the
lantern relation. Hence we have
\[t_{a_r}t_{\omega_k}=t_{h_k}t_\omega t_{\nu_{r,k}}.\]
This proves that $v_k^2=t_{h_k}^{-1}\in K$.
\subsection{Case of $v_kt_{c_{r}}v_k^{-1}$, $g$ -- even}
By the relation
\[v_k^{-1}(c_r)=t_{a_r}t_{e_{r,k-1}}t_{e_{r,k}}^{-1}t_{a_r}^{-1}(c_r), \] we have $v_k^{-1}t_{c_{r}}v_k\in
K$. Since we proved that $v_k^2\in K$, this implies that $v_kt_{c_{r}}v_k^{-1}\in K$.
\subsection{Case of $v_kt_{f_j}v_k^{-1}$}
By the relation $v_k^{-1}(f_j)=t_{f_k}(f_j)$, we have $v_k^{-1}t_{f_j}v_k\in K$ for
$j=1,\ldots,k-1$. Since $v_k^2\in K$, this implies that $v_kt_{f_j}v_k^{-1}\in K$.
\subsection{Case of $v_k v_j v_k^{-1}$}
Using Propositions \ref{pre:prop:slide:cong} and \ref{pre:prop:slide:twist}, it is
straightforward to check that
\[v_j^{-1}(v_k v_j v_k^{-1})=t_{\varepsilon_{k,j}}^{-1},\]
where $\eps_{k,j}$ is the circle indicated in Figure \ref{fig:15_GenCl}(ii) for
$j=k+1,\ldots,n$.
\begin{figure}[h]
\includegraphics{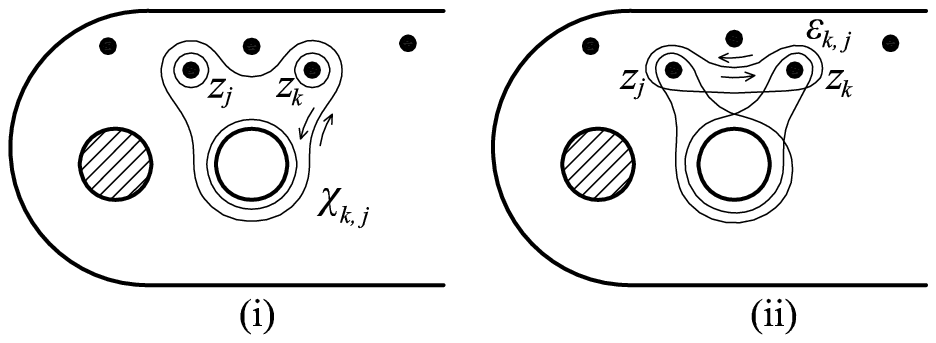}
\caption{Circles of the lantern relation
$t_{a_r}t_{\chi_{k,j}}=t_{\nu_{r,j}}t_{\nu_{r,k}}t_{\eps_{k,j}}$.}\label{fig:15_GenCl}
\end{figure}
Therefore it is enough to prove that $t_{\eps_{k,j}}\in K$. Let
\[\chi_{k,j}=t_{e_{r,k}}^{-1}t_{e_{r,j}}^{-1}t_{e_{r,j-1}}t_{a_r}^{-1}(e_{r,k-1}).\]
Clearly $t_{\chi_{k,j}}\in K$. Now the seven circles indicated in Figure \ref{fig:15_GenCl}
form a configuration of the lantern relation. Hence we have
\[t_{a_r}t_{\chi_{k,j}}=t_{\nu_{r,j}}t_{\nu_{r,k}}t_{\eps_{k,j}}.\]
This proves that $t_{\eps_{k,j}}\in K$.
\subsection{Case of $v_k y v_k^{-1}$}
Let $y'=t_{\nu_{r,k}}v_k^{-1}yv_kt_{\nu_{r,k}}^{-1}$. Since $v_k^2,t_{\nu_{r,k}}\in K$,
it is enough to prove that $y'\in K$. By the relation
$t_{\nu_{r,k}}(v_k^{-1}(\xi))\simeq\xi$, we have $y'(\xi)\simeq\xi$. Therefore we can
assume that $y'(\xi)=\xi$. Let $N_1$ be this of two connected components of $\kre{N\bez \xi}$
which is a support for $y'$ (i.e. $y'$ acts as the identity on the second component --
since $y'$ is conjugate to $y$ this component is well defined). It is known that the
mapping class group ${\cal{M}}(N_1)$ is generated by $y$ and $t_{a_r}$ (cf Theorem A.7 of
\cite{Stukow_twist}), hence $y'$ as a composition of these elements is in $K$.
\end{proof}

\section{Generators for groups ${\cal{PM}}^k(N_{g,s}^n)$ and ${\cal{M}}(N_{g,s}^n)$}
The following proposition is well known and easy to prove, hence we state it without proof.
\begin{prop} \label{prop:gen:ext}
Suppose that we have a short exact sequence of groups and homomorphisms
\[\begin{CD}1@>>>H@>i>> G@>p>> K@>>>1\end{CD}\]
and let $X_H$ and $X_K$ be generating sets for $H$ and $K$ respectively. Let $\fal{X}_K$
be any subset of $G$ such that $p(\fal{X}_K)=X_K$. Then $G$ is generated by
$i(X_H)\cup\fal{X}_K$. \qed
\end{prop}

Let $\fal{N}=N_{g}^{n+s}$ be the surface obtained from $N=N_{g,s}^n$ by gluing a disk
with one puncture to each boundary component and let $1\leq k\leq n$ be an integer. We choose the notation in such a way that
the first $s$ of $n+s$ punctures of $\fal{N}$ correspond to the boundary components of
$N$. Identifying $N$ with a subsurface of $\fal{N}$, we can consider the circles
in ${\cal{C}}$, $\lst{f}{s+k}$ and $\lambda$ as circles on $N$. Similarly we will consider
the puncture slides $v_{s+k+1},\ldots,v_{s+n}$ and the crosscap slide $y$ as elements of
${\cal{PM}}^k(N)$.

By Theorem 2.2 of \cite{Stukow_commen}, we have an exact sequence
\[\begin{CD} 1 @>>> \zz^s @>>> {{\cal{PM}}^k(N)}@>i_{\star}>>{{\cal{PM}}^{s+k}(\fal{N})} @>>>1 \end{CD} \]
where $i_*$ is the homomorphism induced by the inclusion $\map{i}{N}{\fal{N}}$. Moreover,
the generators of $\ker i_*\cong \zz^s$ correspond to the boundary twists
$t_{u_1},\ldots,t_{{u}_{s}}$ on $N$. Let
${\cal{C}}'={\cal{C}}\cup\{\lst{f}{s+k}\}\cup\{\lst{u}{s}\}$. Theorem
\ref{tw:gen:pure:cl} together with Proposition \ref{prop:gen:ext} implies the following.
\begin{tw}\label{tw:gen:pure:bd}
Let $g\geq 3$. Then the mapping class group ${\cal{PM}}^k(N_{g,s}^n)$ is generated by
\begin{itemize}
 \item $\{t_l,v_{s+k+1},\ldots,v_{s+n},y\st{l\in{\cal{C}}'}\}$ if $g$ is odd and
 \item $\{t_l,v_{s+k+1},\ldots,v_{s+n},y,t_{\lambda}\st{l\in{\cal{C}}'}\}$ if $g$ is even.
\end{itemize}\qed
\end{tw}

\begin{figure}[h]
\includegraphics{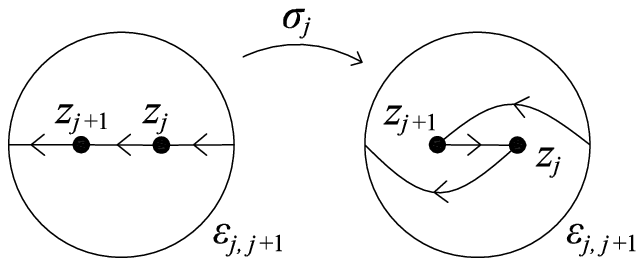}
\caption{Elementary braid $\sig_j$.}\label{fig:16_GenBD}
\end{figure}

Now let us turn to the group ${\cal{M}}(N_{g,s}^n)$. Let $\sig_{s+1},\ldots,\sig_{s+n-1}$
be elementary braids on ${N}$, such that $\sig_j^2=t_{\eps_{j,j+1}}$ (cf Figure
\ref{fig:16_GenBD}), where $\eps_{j,j+1}$ is a circle defined by Figure
\ref{fig:15_GenCl}(ii).
Let ${\cal{C}}'$ be defined as before, i.e.
${\cal{C}}'={\cal{C}}\cup\{\lst{f}{s}\}\cup\{\lst{u}{s}\}$.
\begin{tw}\label{tw:gen:notpure}
Let $g\geq 3$ and $n\geq 2$. Then the mapping class group ${\cal{M}}(N_{g,s}^n)$ is
generated by
\begin{itemize}
 \item $\{t_l,v_{s+1},\sig_{s+1},\ldots,\sig_{s+n-1},y\st{l\in{\cal{C}}'}\}$ if $g$ is
 odd and
 \item $\{t_l,v_{s+1},\sig_{s+1},\ldots,\sig_{s+n-1},y,t_{\lambda}\st{l\in{\cal{C}}'}\}$ if $g$ is even.
\end{itemize}
\end{tw}
\begin{proof}
From the short exact sequence
\[\begin{CD} 1 @>>> {{\cal{PM}}(N_{g,s}^n)} @>i>> {{\cal{M}}(N_{g,s}^n)}@>p>>S_n @>>>1 \end{CD} \]
where $S_n$ is the symmetric group on $n$ letters, by Theorem \ref{tw:gen:pure:bd},
Proposition \ref{prop:gen:ext} and by the fact that
$p(\sig_{s+1}),\ldots,p(\sig_{s+n-1})$ generate $S_n$, we conclude that
${\cal{M}}(N_{g,s}^n)$ is generated by
\begin{itemize}
 \item $\{t_l,v_{s+1},\ldots,v_{s+n},\sig_{s+1},\ldots,\sig_{s+n-1},y\st{l\in{\cal{C}}'}\}$ if $g$ is odd and
 \item $\{t_l,v_{k+1},\ldots,v_{s+n},\sig_{s+1},\ldots,\sig_{s+n-1},y,t_{\lambda}\st{l\in{\cal{C}}'}\}$ if $g$ is even.
\end{itemize}
Since $v_j=\sig_{j-1}^{-1}v_{j-1}\sig_{j-1}$, for $j=s+2,\ldots,s+n$, the generators
$v_{s+2},\ldots,v_{s+n}$ can be removed from the above generating sets.
\end{proof}

\section{Homological results for mapping class groups} The main goal of this section is
to compute the first homology groups of the mapping class groups
${{\cal{PM}}^k(N_{g,s}^n)}$ and ${{\cal{M}}(N_{g,k}^n)}$ for $g\geq 3$. Before we do this, we
need some technical preparations.

Let $N=N_{g,s}^n$,  and $G={{\cal{PM}}^k(N)}$. Moreover, for $f\in G$ let $[f]$ denote the
homology class of $f$ in $\mathrm{H_1}(G)$ (we will use the additive notation in
$\mathrm{H_1}(G)$).
\subsection{Homology classes of twists}
It is well known that two right twists about nonseparating circles on an oriented surface
are conjugate in the mapping class group of this surface. The description of conjugacy
classes of twists about nonseparating circles on a nonorientable surface is slightly more
difficult -- it can be shown \cite{Szep_curv} that there are at least $2^{k+s-1}+1$ such classes in the
group ${{\cal{PM}}^k(N_{g,s}^n)}$. However, for our purposes it will
suffice to consider a very simple case of this description, namely when the complement of
a circle in $N$ is nonorientable -- see Proposition \ref{prop:tw:nonsep:cong}.

\begin{lem}\label{lem:one:orbit:circl}
Let $c$ be a nonseparating two--sided circle on $N$ such that $N\bez c$ is nonorientable.
Then there exists $f\in G$ such that $f(c)=a_1$ or $f(c)=a_1^{-1}$, where $a_1$ is as in Figures
\ref{fig:02_GenRed} and \ref{fig:03_GenRed}.
\end{lem}
\begin{proof}
Clearly $N\bez c$ and $N\bez a_1$ are diffeomorphic. It is a standard argument that we
can choose this diffeomorphism in such a way that it extends to a diffeomorphism
$\map{f}{N}{N}$ such that $f\in G$ and $f(c)=a_1^{\pm1}$.
\end{proof}

\begin{lem}\label{lem:cong:bu:kl}
Let $c$ be a nonseparating two--sided circle on a Klein bottle $N=N_{2,1}$ with one
boundary component. Then there exists $f\in {{\cal{M}}(N)}$ such that
$ft_cf^{-1}=t_c^{-1}$.
\end{lem}
\begin{proof}
The assertion follows from the structure of ${{\cal{M}}(N_{2,1})}$ -- cf Theorem A.7 of
\cite{Stukow_twist}.
\end{proof}
\begin{lem}\label{lem:cong:nonsep}
Let $c$ be a nonseparating two--sided circle on $N$ such that $N\bez c$ is nonorientable.
Then $t_c$ and $t_{c}^{-1}$ are conjugate in $G$. In particular $2[t_c]=[{t}_c^2]=0$.
\end{lem}
\begin{proof}
By Lemma \ref{lem:one:orbit:circl}, $t_c$ is conjugate to $t_{a_1}$ or to $t_{a_1}^{- 1}$. Hence the
assertion follows form Lemma \ref{lem:cong:bu:kl} applied to the Klein bottle cut off by
a circle $\xi$ in Figure \ref{fig:05_GenRed}.
\end{proof}
\begin{prop} \label{prop:tw:nonsep:cong}
Let $c$ be a nonseparating two--sided circle on $N$ such that $N\bez c$ is nonorientable.
Then the twist $t_{c}$ is conjugate to $t_{a_1}$ in $G$. In particular $[t_c]=[t_{a_1}]$
in $\mathrm{H_1}(G)$.
\end{prop}
\begin{proof}
By Lemma \ref{lem:one:orbit:circl}, $t_{c}$ is conjugate to $t_{a_1}$ or to
$t_{a_1}^{-1}$. Hence Lemma~\ref{lem:cong:nonsep} implies that $t_{c}$ is conjugate to
$t_{a_1}$.
\end{proof}
\begin{lem}\label{lem:cong:br}
Assume $g=2r+2\geq 4$. Then $t_{b_{r+1}}$ and $t_{b_{r+1}}^{-1}$ are conjugate in $G$,
where $b_{r+1}$ is as in Figure \ref{fig:03_GenRed}. In particular
$2[t_{b_{r+1}}]=[t_{b_{r+1}}^2]=0$.
\end{lem}
\begin{proof}
To obtain the conclusion it is enough to apply Lemma \ref{lem:cong:bu:kl} to the Klein
bottle cut off by a circle $\xi'$ indicated in Figure \ref{fig:17_Homo}.
\begin{figure}[h]
\includegraphics{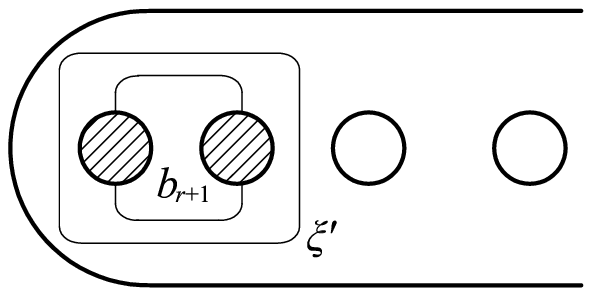}
\caption{Circles ${b_{r+1}}$ and ${\xi'}$ -- Lemma \ref{lem:cong:br}.}\label{fig:17_Homo}
\end{figure}
\end{proof}
\begin{lem}\label{lem:cong:b_r:a_r}
Assume $g=2r+2\geq 6$. Then $[t_{b_{r+1}}]=0$.
\end{lem}
\begin{proof}
Figure \ref{fig:18_Homo} shows that there is a lantern configuration with one twist
$t_{b_{r+1}}$ and each of the remaining six circles is conjugate to $t_{a_1}$ (cf Proposition
\ref{prop:tw:nonsep:cong}).
\begin{figure}[h]
\includegraphics{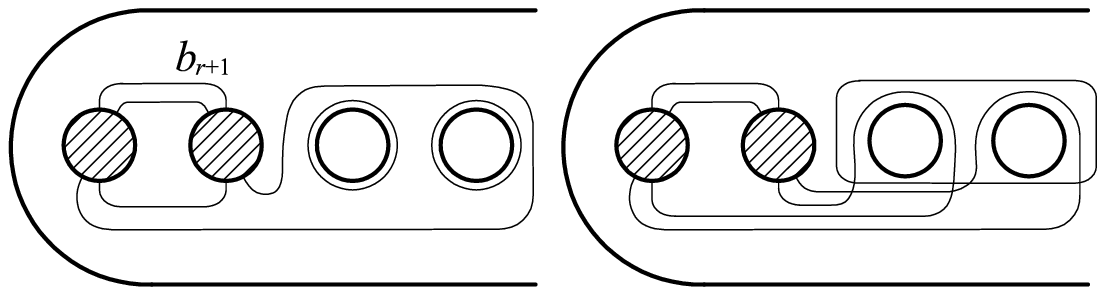}
\caption{Lantern relation -- Lemma \ref{lem:cong:b_r:a_r}.}\label{fig:18_Homo}
\end{figure}
Hence we have the relation (in $\mathrm{H_1}(G)$)
\[ [{t_{b_{r+1}}}{t_{a_1}}{t_{a_1}}{t_{a_1}}]=[{t_{a_1}}{t_{a_1}}{t_{a_1}}]. \]
\end{proof}
\begin{lem}\label{lem:cong:twists:7}
Assume $g\geq 7$. Then $[t_{a_{1}}]=0$.
\end{lem}
\begin{proof}
Figure \ref{fig:19_Homo} shows that there is a lantern configuration with all twists
conjugate to $t_{a_1}$.
\begin{figure}[h]
\includegraphics{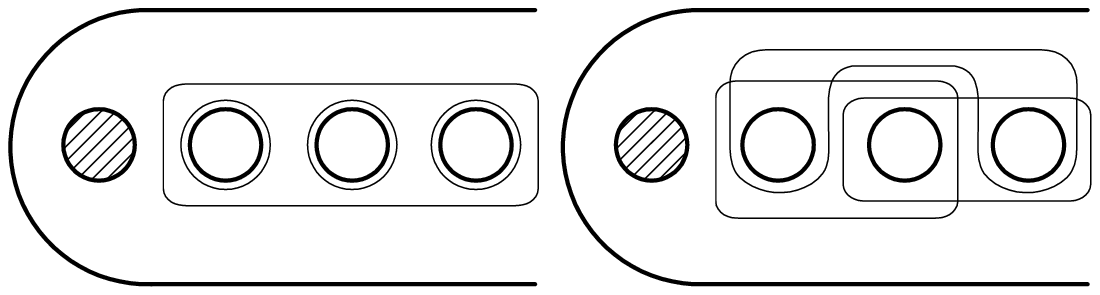}
\caption{Lantern relation -- Lemma \ref{lem:cong:twists:7}.}\label{fig:19_Homo}
\end{figure}
Hence we have the relation
\[ [{t_{a_{1}}}{t_{a_1}}{t_{a_1}}{t_{a_1}}]=[{t_{a_1}}{t_{a_1}}{t_{a_1}}]. \]
\end{proof}

\subsection{Homology classes of crosscap slides}
\begin{lem}\label{lem:cong:kl:y}
Let $y$ be a crosscap slide on a Klein bottle $N=N_{2,1}$ with one boundary component.
Then there exists a diffeomorphism $\map{f}{N}{N}$ such that $fyf^{-1}=y^{-1}$ and
$f|_{\partial N}=-id$.
\end{lem}
\begin{proof}
The lemma can be easily deduced from the proof of Theorem 5.8 of \cite{Kork-non}. It is
also a direct consequence of the structure of ${\cal{M}}(N)$ -- cf Theorem A.7 of
\cite{Stukow_twist}.
\end{proof}
\begin{lem}\label{lem:cong:crosscap}
Let $y$ be a crosscap slide on $N$ such that $y^2$ is a twist $t_\xi$ about a circle
which separates $N$ into two nonorientable components. Then $y$ and $y^{-1}$ are
conjugate in $G$. In particular $2[y]=[y^2]=0$.
\end{lem}
\begin{proof}
Let $N_1$ and $N_2$ be components of $\kre{N\bez \xi}$ and assume that $N_1$ is the
support for $y$. By Lemma \ref{lem:cong:kl:y}, it is enough to prove that there exists a
diffeomorphism of $N_2$ which preserves the local orientation around each of the
punctures, acts as the $-id$ on $\partial N_2$ coming from $\xi$ and as the $id$ on each
of the remaining boundary components. Since $N_2$ is nonorientable such diffeomorphism
can be obtained by sliding the boundary component coming from $\xi$ along a one--sided
loop.
\end{proof}
\subsection{Homology classes of puncture slides}
\begin{lem}\label{lem:sqr:punc:slide}
Assume that $g\geq 3$. Then $2[v_j]=[v_j^2]=[t_{h_j}^{-1}]=0$, for $j=s+k+1,\ldots,s+n$.
\end{lem}
\begin{proof}
Let us recall the instance of the lantern relation from Section~\ref{sec:v_k} (cf Figure
\ref{fig:14_GenCl}). After replacing $k$ with $j$ we can write it as
\[t_{a_r}t_{\omega_j}=t_{h_j}t_\omega t_{\nu_{r,j}}.\]
By Proposition \ref{prop:tw:nonsep:cong},
$[t_{a_r}]=[t_{\omega_j}]=[t_\omega]=[t_{\nu_{r,j}}]$. Hence $[t_{h_j}]=0$.
\end{proof}
\begin{lem}\label{lem:cong:punct:slides}
Assume that $n\geq 2$. Then all the puncture slides $v_{s+1},\ldots,v_{s+n}$ are contained in the same conjugacy class in ${\cal{M}}(N)$. In particular all these elements are equal in
$\mathrm{H_1}({\cal{M}}(N))$.
\end{lem}
\begin{proof}
The assertion follows inductively by the relation
\[ v_{j+1}=\sig_j^{-1} v_j\sig_j\quad\text{for $j=s+1,\ldots,s+n-1$}.\]
\end{proof}
\subsection{Homology classes of boundary twists}
\begin{lem}\label{lem:triv:bound:tw}
Assume that $g\geq 5$. Then the boundary twists $t_{u_1},\ldots,t_{u_s}$ are trivial in
$\mathrm{H_1}(G)$.
\end{lem}
\begin{proof}
Figure \ref{fig:20_Homo} shows that for each $1\leq j\leq s$, there exists a lantern
configuration with one twist $t_{u_j}$ and all other twists conjugate to $t_{a_1}$.
\begin{figure}[h]
\includegraphics{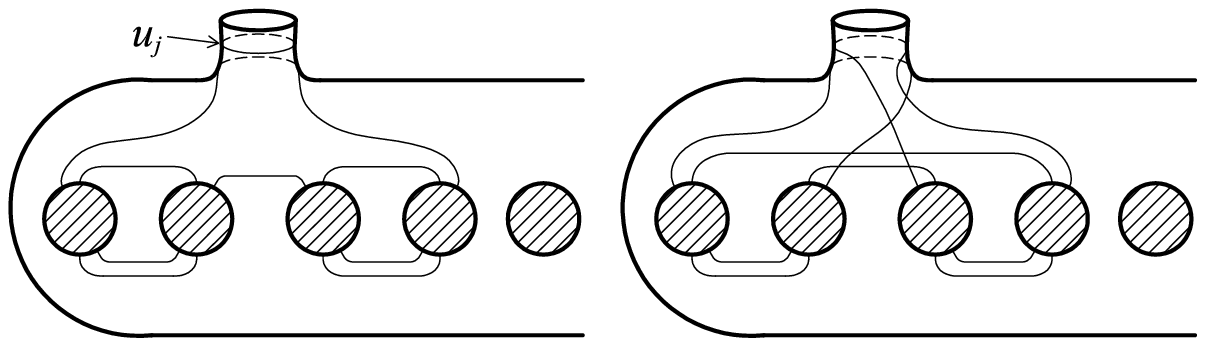}
\caption{Lantern relation -- Lemma \ref{lem:triv:bound:tw}.}\label{fig:20_Homo}
\end{figure}
Hence we have
\[
[t_{u_j}t_{a_1}t_{a_1}t_{a_{1}}]=[t_{a_1}t_{a_1}t_{a_{1}}].\]
\end{proof}
\begin{lem}\label{lem:triv:kap}
Let $\kappa$ be a circle on $N$ as in Figure \ref{fig:21_Homo} for $g=3$ and as in Figure
\ref{fig:21_2_Homo} for $g=4$. Then $[t_{\kappa}]=0$.
\begin{figure}[h]
\includegraphics{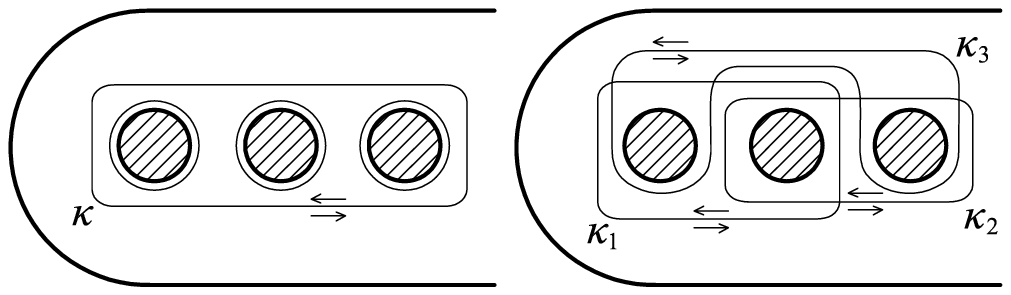}
\caption{Circles of the lantern relation
$[t_\kappa]=[t_{\kappa_1}t_{\kappa_2}t_{\kappa_3}]$.}\label{fig:21_Homo}
\end{figure}
\begin{figure}[h]
\includegraphics{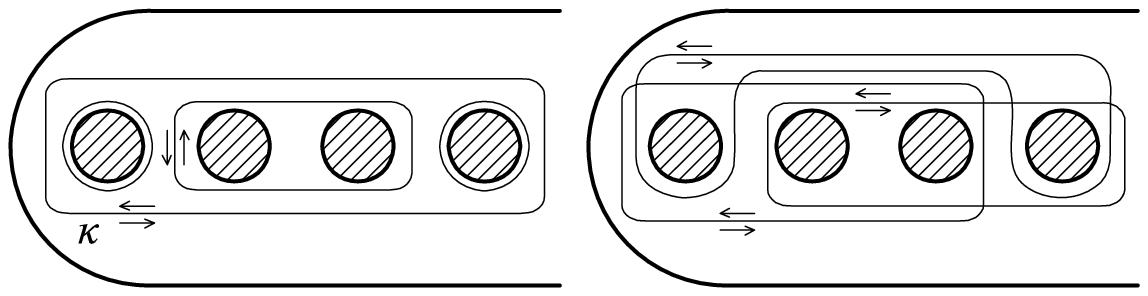}
\caption{Lantern relation for $g=4$, Lemma \ref{lem:triv:kap}.}\label{fig:21_2_Homo}
\end{figure}
\end{lem}
\begin{proof}
Suppose first that $g=3$ and let the circles $\kappa_1,\kappa_2,\kappa_3$ be as in Figure
\ref{fig:21_Homo}. For each $i=1,2,3$ we have $t_{\kappa_i}=y_i^2$, where $y_i$ is a
crosscap slide satisfying the assumptions of Lemma \ref{lem:cong:crosscap}. In fact, we can define $y_i$ to be a crosscap slide on the component of $N\bez \kappa_i$ which is diffeomorphic to a Klein bottle with one boundary component. By Lemma \ref{lem:cong:crosscap}, $[t_{\kappa_1}]=[t_{\kappa_2}]=[t_{\kappa_3}]=0$. Moreover, Figure
\ref{fig:21_Homo} shows that there is a lantern relation
\[[t_\kappa]=[t_{\kappa_1}t_{\kappa_2}t_{\kappa_3}].\]
The case $g=4$ is very similar. Observe that by the reasoning for $g=3$, the homology
classes of all the circles but $\kappa$ in Figure \ref{fig:21_2_Homo} are trivial. Hence
the homology class of $\kappa$ is also trivial.
\end{proof}

By a \emph{region} of a surface $N$ we will mean any closed connected subsurface $\Delta$ of $N$ of
genus $0$ which satisfies the following two conditions.
\begin{enumerate}
 \item $\Delta$ has one distinguished boundary component which is disjoint from the set of punctures and from the boundary of $N$. We will denote this distinguished boundary component by $\partial \Delta$.
 \item Every boundary component of $\Delta$ different from $\partial \Delta$ is a boundary component of $N$.
\end{enumerate}
In other words $\Delta$ is a disk with punctures and/or boundary components of $N$ imbedded in
$N$.
\begin{lem}\label{lem:disk:im:sur}
Let $\Delta$ be a region of a surface $N=N_{g,s}^n$ for $g\geq 3$, such that
\begin{itemize}
 \item $\KL{j\st j\in \{1,\ldots,s\} \text{ and } u_j\subset\Delta}=\Theta$,
 \item $\KL{j\st j\in \{s+1,\ldots,s+n\}\text{ and } z_j\in\Delta}=\Omega$,
 \item if $\Theta\neq\emptyset$ then the orientation of $\Delta$ agrees with the orientation of a neighbourhood of
 $u_\theta$ for $\theta\in \Theta$ (in other words $t_{u_\theta}$ is a right twist on~$\Delta$).
 \item the orientation of a neighbourhood of $\partial \Delta$ agrees with the orientation of $\Delta$.
\end{itemize}
Then
\[[t_{\partial \Delta}]=\left[\prod_{\theta \in \Theta}t_\theta \right]=\sum_{\theta \in \Theta}[t_\theta].\]
\end{lem}
\begin{proof}
The proof is by induction on $|\Omega\cup\Theta|$. For $|\Omega\cup\Theta|=1$ there is
nothing to prove, and to make the inductive step, assume that $\Delta'$ is obtained from
$\Delta$ by adding one puncture/boundary component. Figure~\ref{fig:22_Homo} shows that
we have the lantern relation
\[[t_{\partial\Delta}t_u t_{a_1}t_{a_1}]=[t_{\partial\Delta'}t_{a_1}t_{a_1}],\]
\begin{figure}[h]
\includegraphics{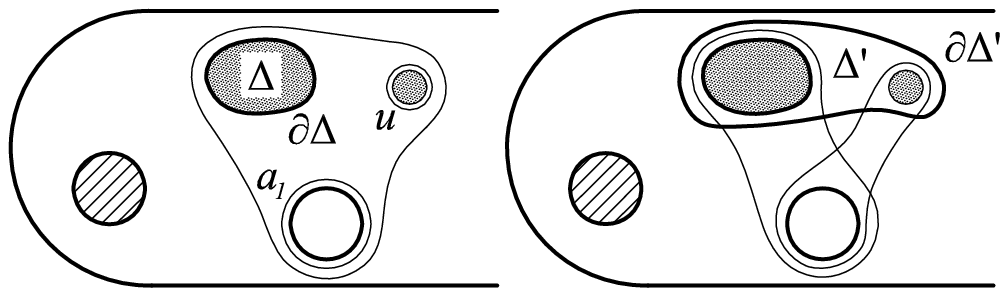}
\caption{Lantern relation $[t_{\partial\Delta}t_u
t_{a_1}t_{a_1}]=[t_{\partial\Delta'}t_{a_1}t_{a_1}]$.}\label{fig:22_Homo}
\end{figure}
where $u$ is a circle around the added puncture/boundary component. Therefore
$[t_{\partial\Delta'}]=[t_{\partial\Delta}]+[t_u]$ and $t_u$ is trivial if we added a puncture.
\end{proof}
\begin{lem} \label{lem:pro:bd:tw:3}
Let $g=3$ or $g=4$, then $[t_{u_1}]+[t_{u_2}]+\cdots +[t_{u_s}]=0$.
\end{lem}
\begin{proof}
By Lemma \ref{lem:disk:im:sur}, the homology class of $t_{u_1}t_{u_2}\cdots t_{u_s}$ is
equal to the homology class of a twist about a circle which is the boundary of a region
containing all the punctures and boundary components, i.e. along the circle $\kappa$ in
Figure \ref{fig:21_Homo} for $g=3$ and Figure \ref{fig:21_2_Homo} for $g=4$ respectively.
By Lemma \ref{lem:triv:kap}, $[t_\kappa]=0$.
\end{proof}
\begin{prop}\label{prop:separ:non:cong}
Let $c$ be a two--sided circle on $N$ such that $c$ separates $N$ into two nonorientable
surfaces. Then $t_c$ and $t_c^{-1}$ are conjugate in $G$. In particular $2[t_c]=[t_c^2]=0$.
\end{prop}
\begin{proof}
Let $N_1$ and $N_2$ be components of $\kre{N\bez c}$. For each $i=1,2$ there exists a
diffeomorphism $f_i$ of $N_i$ which preserves the local orientation around each of the
punctures, acts as the $-id$ on $\partial N_i$ coming from $c$ and as the $id$ on each of
the remaining boundary components. In fact, since $N_i$ is nonorientable such a
homeomorphism can be obtained by sliding the boundary component coming from $c$ along a
one--sided loop. Now if $f$ is a diffeomorphism of $N$ obtained by gluing $f_1$ and
$f_2$, then $ft_cf^{-1}=t_c^{-1}$.
\end{proof}
\begin{lem} \label{lem:sq:tr:boun:tws}
Let $g\geq 3$, then $2[t_{u_j}]=[t_{u_j}^2]=0$.
\end{lem}
\begin{proof}
Figure \ref{fig:23_Homo} shows that there is a lantern relation
\[[t_{u_j}t_{a_1}t_{a_1}]=[t_{\eta_j}t_{a_1}t_{a_1}].\]
\begin{figure}[h]
\includegraphics{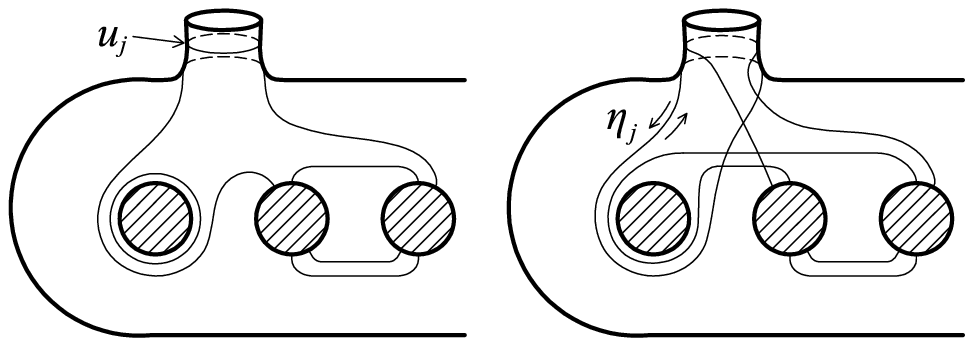}
\caption{Lantern relation
$[t_{u_j}t_{a_1}t_{a_1}]=[t_{\eta_j}t_{a_1}t_{a_1}]$.}\label{fig:23_Homo}
\end{figure}
Hence $[t_{u_j}]=[t_{\eta_j}]$ and by Proposition \ref{prop:separ:non:cong},
$[t_{\eta_j}^2]=0$.
\end{proof}
\begin{uw}\label{rem:same:notpure}
Observe that since the identity $\map{i}{N_{g,s}^n}{N_{g,s}^n}$ induces a homomorphism
$\map{i_*}{{\cal{PM}}^k(N_{g,s}^n)}{{\cal{M}}(N_{g,s}^n)}$, Lemmas \ref{lem:cong:nonsep},
\ref{lem:cong:br}, \ref{lem:cong:b_r:a_r}, \ref{lem:cong:twists:7},
\ref{lem:cong:crosscap}, \ref{lem:sqr:punc:slide}, \ref{lem:triv:bound:tw},
\ref{lem:disk:im:sur}, \ref{lem:pro:bd:tw:3}, \ref{lem:sq:tr:boun:tws} and Proposition
\ref{prop:tw:nonsep:cong} remain true if we replace the group ${{\cal{PM}}^k(N_{g,s}^n)}$
with ${{\cal{M}}(N_{g,s}^n)}$.
\end{uw}
\subsection{Homology classes of elementary braids}
\begin{lem} \label{lem:braid:cong}
Let $\sig_{s+1},\ldots,\sig_{s+n-1}$ be elementary braids on $N=N_{g,s}^n$ as in Theorem
\ref{tw:gen:notpure}. Then $\sig_{s+1},\ldots,\sig_{s+n-1}$ are in a single conjugacy
class in ${\cal{M}}(N)$. In particular all these elements are equal in
$\mathrm{H_1}({\cal{M}}(N))$.
\end{lem}
\begin{proof}
The assertion follows inductively from the braid relation
\[\sig_{j+1}=(\sig_j\sig_{j+1})\sig_j(\sig_j\sig_{j+1})^{-1}\quad\text{for
$j=s+1,\ldots,s+n-2$}.\]
\end{proof}
\begin{lem}\label{lem:braid:sq:trv}
Let $g\geq 3$. Then $[\sig_{s+1}^2]=0$ in $\mathrm{H_1}({\cal{M}}(N))$.
\end{lem}
\begin{proof}
Since $\sig_{s+1}^2$ is a twist about the boundary of a region containing the punctures
$z_{s+1}$ and $z_{s+2}$, the assertion follows from Lemma \ref{lem:disk:im:sur} (cf
Remark \ref{rem:same:notpure}).
\end{proof}
\subsection{Main theorems}
\begin{tw}\label{tw:hom:pure}
Let $N=N_{g,s}^n$. Then \[\mathrm{H_1}({\cal{PM}}^k(N),\zz)
\cong\begin{cases}\zz_2^{2+n-k}&\text{if $g=3$ and $s=0$},\\
\zz_2^{1+n-k+s}&\text{if $g=3$ and $s\geq 1$},\\
\zz_2^{3+n-k}&\text{if $g=4$ and $s=0$},\\
\zz_2^{2+n-k+s}&\text{if $g=4$ and $s\geq 1$},\\
\zz_2^{2+n-k}&\text{if $g=5,6$},\\
\zz_2^{1+n-k}&\text{if $g\geq 7$}.\\
\end{cases}\]
\end{tw}
\begin{proof}
By Theorem \ref{tw:gen:pure:bd} and Proposition \ref{prop:tw:nonsep:cong},
$\mathrm{H_1}({\cal{PM}}^k(N))$ is generated by
\begin{itemize}
 \item $[t_{a_1}],[y],[v_{s+k+1}],\ldots,[v_{s+n}],[t_{u_1}],\ldots,[t_{u_{s}}]$ if $g$ is
 odd and
 \item $[t_{a_1}],[t_{b_{r+1}}],[y],[v_{s+k+1}],\ldots,[v_{s+n}],[t_{u_1}],\ldots,[t_{u_{s}}]$ if
 $g$ is even.
\end{itemize}
By Lemmas \ref{lem:cong:nonsep}, \ref{lem:cong:br}, \ref{lem:cong:crosscap},
\ref{lem:sqr:punc:slide} and \ref{lem:sq:tr:boun:tws}, each of these generators has order
at most $2$. Moreover,
\begin{itemize}
 \item $[t_{a_1}]=0$ if $g\geq 7$ (Lemma \ref{lem:cong:twists:7}),
 \item $[t_{b_{r+1}}]=0$ if $g\geq 6$ (Lemma \ref{lem:cong:b_r:a_r}),
 \item $[t_{u_1}]=\ldots=[t_{u_s}]=0$ if $g\geq 5$ (Lemma \ref{lem:triv:bound:tw}),
 \item $[t_{u_s}]=-([t_{u_1}]+\cdots +[t_{u_{s-1}}])$ if $g=3,4$ (Lemma \ref{lem:pro:bd:tw:3}).
\end{itemize} Therefore $\mathrm{H_1}({\cal{PM}}^k(N))$ is generated by
\begin{itemize}
 \item $[t_{a_1}],[y],[v_{s+k+1}],\ldots,[v_{s+n}],[t_{u_1}],\ldots,[t_{u_{s-1}}]$ if $g=3$,
 \item $[t_{a_1}],[t_{b_{r+1}}],[y],[v_{s+k+1}],\ldots,[v_{s+n}],[t_{u_1}],\ldots,[t_{u_{s-1}}]$ if $g=4$,
 \item $[t_{a_1}],[y],[v_{s+k+1}],\ldots,[v_{s+n}]$ if $g=5,6$,
 \item $[y],[v_{s+k+1}],\ldots,[v_{s+n}]$ if $g\geq 7$.
\end{itemize}
To finish the proof it is enough to show that these elements are independent over
$\zz_2$.

Suppose that
\begin{equation}
[t_{a_1}^{\alpha}t_{b_{r+1}}^\beta y^\gamma v_{s+k+1}^{\nu_{s+k+1}}\cdots
v_{s+n}^{\nu_{s+n}}t_{u_1}^{\mu_1}\cdots t_{u_{s-1}}^{\mu_{s-1}}]=0,
\label{eq:triv}
\end{equation}
for some
$\alpha,\beta,\gamma,\nu_i,\mu_l\in\zz_2$ such that $\alpha=0$ for $g\geq 7$, $\beta=0$ for $g\neq 4$ and $\mu_{1}=\ldots=\mu_{s-1}=0$ for $g\geq 5$. Our first goal is to show that $\alpha=\beta=\gamma=0$.

Let $\fal{N}$ be the closed surface obtained from $N$ by forgetting the punctures and
gluing a disk to each boundary component of $N$. The inclusion $N\hookrightarrow \fal{N}$ induces the homomorphism
\[\map{\Phi}{{\cal{M}}(N)}{{\cal{M}}(\fal{N})}.\]
Clearly $\Phi(v_i)=\Phi(t_{u_l})=1$, hence equation \eqref{eq:triv} yields
\[ \alpha [t_{a_1}]+\beta [t_{b_{r+1}}]+\gamma [y]=[t_{a_1}^{\alpha}t_{b_{r+1}}^\beta y^\gamma]=0\quad \text{in $\mathrm{H_1}({{\cal{M}}(\fal{N})})$}.\]
By Theorem 1.1 of \cite{Kork-non1},
\[\mathrm{H_1}({\cal{M}}(\fal{N}))\cong\begin{cases}
\gen{t_{a_1},y}\cong\zz_2\oplus\zz_2&\text{if $g=3,5,6$},\\
\gen{t_{a_1},t_{b_{r+1}},y}\cong\zz_2\oplus\zz_2\oplus\zz_2&\text{if $g=4$},\\
\gen{y}\cong\zz_2&\text{if $g\geq 7$}.
\end{cases}
\]
Hence $\alpha=\beta=\gamma=0$ and equation \eqref{eq:triv} takes form
\begin{equation}\label{eq:triv2}
[v_{s+k+1}^{\nu_{s+k+1}}\cdots v_{s+n}^{\nu_{s+n}}t_{u_1}^{\mu_1}\cdots
t_{u_{s-1}}^{\mu_{s-1}}]=0.
\end{equation}

In order to show that $\nu_{s+k+1}=\cdots=\nu_{s+n}=0$, define \[\map{\Psi_j}{{\cal{M}}(N)}{\zz_2},\] for $j=s+k+1,\ldots,s+n$, as follows:
$\Psi_j(f)=1$ if $f$ reverses the local orientation around the puncture $z_j$, and
$\Psi_j(f)=0$ otherwise. Since $\Psi_j(v_i)=0$ for $i\neq j$ and $\Psi_j(t_{u_l})=0$,
equation \eqref{eq:triv2} implies that
\[\nu_j=\nu_j[\Psi_j(v_j)]=[\Psi_j(v_j^{\nu_j})]=0 \quad \text{in $\mathrm{H_1}(\zz_2)=\zz_2$}.\]
Therefore equation \eqref{eq:triv2} becomes
\begin{equation}\label{eq:triv3}
[t_{u_1}^{\mu_1}\cdots t_{u_{s-1}}^{\mu_{s-1}}]=0.
\end{equation}
This completes the proof for $g\geq 5$, hence assume that $g=3$ or $g=4$. We can also assume that $s\geq 2$.

Let $\dasz{N}_j$, for $j=1,\ldots,{s-1}$, be the surface obtained from $N$ by
forgetting the punctures, gluing a cylinder to the circles $u_j$ and $u_s$ and finally, gluing a disk to each of the remaining boundary components. Then $\dasz{N}_j$ is a nonorientable surface of genus either $5$ or $6$, with neither punctures nor boundary components.
Let
\[\map{\Upsilon_j}{{\cal{M}}(N)}{{\cal{M}}(\dasz{N}_j)}\]
be the homomorphism induced by inclusion. Since $\Upsilon_j(t_{u_l})=0$ for $l\neq j$
and $l\neq s$, equation \eqref{eq:triv3} gives us
\[\mu_j [t_{u_j}]=[t_{u_j}^{\mu_j}] =0\quad \text{in $\mathrm{H_1}({{\cal{M}}(\dasz{N}_j)})$}.\]
Since ${u_j}$ is a nonseparating circle on $\dasz{N}_j$ and $\dasz{N}_j\bez {u_j}$ is nonorientable, by Theorem 1.1 of \cite{Kork-non1}, $[t_{u_j}]\neq 0$ in
$\mathrm{H_1}({{\cal{M}}(\dasz{N}_j)})$, hence $\mu_j=0$.
\end{proof}

\begin{tw}\label{tw:hom:notpure}
Let $N=N_{g,s}^n$ with $n\geq 2$. Then
\[\mathrm{H_1}({\cal{M}}(N),\zz)
\cong\begin{cases}\zz_2^4&\text{if $g=3$ and $s=0$},\\
\zz_2^{s+3}&\text{if $g=3$ and $s\geq 1$},\\
\zz_2^{5}&\text{if $g=4$ and $s=0$},\\
\zz_2^{s+4}&\text{if $g=4$ and $s\geq 1$},\\
\zz_2^{4}&\text{if $g=5,6$},\\
\zz_2^{3}&\text{if $g\geq 7$}.\\
\end{cases}\]
\end{tw}
\begin{proof}
The proof follows similar lines to the proof of Theorem \ref{tw:hom:pure}. By Theorem
\ref{tw:gen:notpure}, Propositions \ref{prop:tw:nonsep:cong} and \ref{lem:braid:cong},
Lemmas \ref{lem:cong:b_r:a_r}, \ref{lem:cong:twists:7}, \ref{lem:cong:punct:slides},
\ref{lem:triv:bound:tw} and \ref{lem:pro:bd:tw:3} (cf Remark \ref{rem:same:notpure}),
$\mathrm{H_1}({\cal{M}}(N))$ is generated by
\begin{itemize}
 \item $[t_{a_1}],[v_{s+1}],[y],[\sig_{s+1}],[t_{u_1}],\ldots,[t_{u_{s-1}}]$ if $g=3$,
 \item $[t_{a_1}],[v_{s+1}],[t_{b_{r+1}}],[y],[\sig_{s+1}],[t_{u_1}],\ldots,[t_{u_{s-1}}]$ if $g=4$,
 \item $[t_{a_1}],[v_{s+1}],[y],[\sig_{s+1}]$ if $g=5,6$,
 \item $[v_{s+1}],[y],[\sig_{s+1}]$ if $g\geq 7$.
\end{itemize}
Moreover, by Lemmas \ref{lem:cong:nonsep}, \ref{lem:cong:br}, \ref{lem:cong:crosscap},
\ref{lem:sqr:punc:slide}, \ref{lem:sq:tr:boun:tws} and \ref{lem:braid:sq:trv} (cf Remark
\ref{rem:same:notpure}), each of these generators has order at most $2$.

As in the proof of Theorem \ref{tw:hom:pure} suppose that
\begin{equation}
[t_{a_1}^{\alpha}t_{b_{r+1}}^\beta y^\gamma v_{s+1}^\nu \sig_{s+1}^{\eps}t_{u_1}^{\mu_1}\cdots
t_{u_{s-1}}^{\mu_{s-1}}]=0, \label{eq:tw2:triv}
\end{equation}
for some $\alpha,\beta,\gamma,\eps,\nu,\mu_i\in\zz_2$. As in the proof
of Theorem \ref{tw:hom:pure} we conclude that $\alpha=\beta=\gamma=\mu_1=\ldots=\mu_{s-1}=0$
and equation \eqref{eq:tw2:triv} becomes
\begin{equation}
[v_{s+1}^\nu \sig_{s+1}^{\eps}]=0. \label{eq:tw2:triv2}
\end{equation}
Let
\[\map{\Theta}{{\cal{M}}(N)}{\zz_2} \]
be the homomorphism defined as follows: $\Theta(f)$ is the sign of a permutation
\[(z_{s+1},\ldots,z_{s+n})\mapsto (f(z_{s+1}),\ldots,f(z_{s+n}))\]
(we use the additive notation for $\zz_2$, hence $0$ and $1$ means even and odd
permutation respectively). Clearly $\Theta(v_{s+1})=0$ and $\Theta(\sig_{s+1})=1$, hence
by equation \eqref{eq:tw2:triv2},
\[\eps=\eps[\Theta(\sig_{s+1})]=[\Theta(\sig_{s+1}^{\eps})]=0 \quad \text{in $\mathrm{H_1}(\zz_2)=\zz_2$}.\]
It remains to show that $\nu=0$. Let $\Delta_{s+1},\ldots,\Delta_{s+n}$ be the
collection of disjoint oriented disks on $N$ such that $\Delta_j\cap
\{z_{s+1},\ldots,z_{s+n}\}=z_j$. Up to isotopy, every $f\in {\cal{M}}(N)$ restricts to
$\bigcup_{j=s+1}^{s+n}\Delta_j$. For each $f\in {\cal{M}}(N)$ define
\[\map{\Psi_j}{{\cal{M}}(N)}{\zz_2}\]
as follows: $\Psi_j(f)=0$ if $\map{f_{|\Delta_j}}{\Delta_j}{f(\Delta_j)}$ is orientation
preserving and $\Psi_j(f)=1$ otherwise (observe that $\Psi_j$ is just a map -- it is not
a homomorphism). Let $\Psi=\Psi_{s+1}+\cdots+\Psi_{s+n}$. It is not difficult to check
that
\[\map{\Psi}{{\cal{M}}(N)}{\zz_2}\]
is a homomorphism. Since $\Psi(v_{s+1})=1$, $[v_{s+1}]\neq 0$ in
$\mathrm{H_1}({\cal{M}}(N))$, hence $\nu=0$.
\end{proof}

\bibliographystyle{abbrv}
\bibliography{mybib}
\end{document}